\RequirePackage{ifpdf}
\ifpdf 
\documentclass[pdftex]{sigma}
\else
\documentclass{sigma}
\fi

\renewcommand{\a }{\alpha }

\renewcommand{\d}{\delta }
\newcommand{\D }{\Delta }

\newcommand{\e }{\varepsilon }
\newcommand{\g }{\gamma}

\renewcommand{\l }{\lambda }
\renewcommand{\L }{\Lambda }

\newcommand{\n }{\nabla }
\newcommand{\var }{\varphi }

\newcommand{\s }{\sigma }
\newcommand{\Sig }{\Sigma}

\newcommand{\ov}{\overline}
\newcommand{\intbar}{\mathop{\int\makebox(-13.5,0){\rule[4pt]{.7em}{0.3pt}}%
\kern-6pt}\nolimits}

\newcommand{\be}{\begin{gather}}
\newcommand{\ee}{\end{gather}}

\newcommand{\R}{\mathbb{R}}

\newcommand{\N}{\mathbb{N}}

\begin{document}

\allowdisplaybreaks
	
\renewcommand{\PaperNumber}{120}

\FirstPageHeading

\renewcommand{\thefootnote}{$\star$}

\ShortArticleName{Conformal Metrics with Constant $Q$-Curvature}

\ArticleName{Conformal Metrics with Constant $\boldsymbol{Q}$-Curvature\footnote{This paper is a
contribution to the Proceedings of the 2007 Midwest
Geometry Conference in honor of Thomas~P.\ Branson. The full collection is available at
\href{http://www.emis.de/journals/SIGMA/MGC2007.html}{http://www.emis.de/journals/SIGMA/MGC2007.html}}}

\Author{Andrea MALCHIODI}

\AuthorNameForHeading{A. Malchiodi}

\Address{SISSA, Via Beirut 2-4, Trieste, Italy}

\Email{\href{mailto:malchiod@sissa.it}{malchiod@sissa.it}}

\URLaddress{\url{http://people.sissa.it/~malchiod/}}

\ArticleDates{Received September 02, 2007, in f\/inal form December
05, 2007; Published online December 13, 2007}

\Abstract{We consider the problem of varying conformally the metric
of a four dimensional manifold in order to obtain constant
$Q$-curvature. The problem is variational, and solutions are in
general found as critical points of saddle type. We show how the
problem leads naturally to consider the set of formal barycenters of
the manifold.}

\Keywords{$Q$-curvature; geometric PDEs; variational methods;
min-max schemes}

\Classification{35B33; 35J35; 53A30; 53C21}

\rightline{\em Dedicated to Tom Branson}

\section{Introduction}\label{s:in}

$Q$-curvature is a scalar quantity introduced by Tom Branson in
1985, \cite{br2}, which in recent years has attracted lot of
attention since it turns out to be fundamental in several questions.
Here we will be mainly interested in its role in conformal geometry,
and especially in its transformation law under a conformal change of
metric.

Let $(M,g)$ be a four dimensional manifold, with Ricci tensor
${\rm Ric}_g$, scalar curvature $R_g$ and Laplace--Beltrami operator
$\D_g$: the $Q$-curvature of $M$ is def\/ined by
\begin{gather*}
    Q_g = - \frac{1}{12} \left( \D_g R_g - R_g^2 + 3 |{\rm Ric}_g|^2
    \right).
\end{gather*}
If we perform the conformal change of metric $\tilde{g} = e^{2w} g$,
for some smooth function $w$ on $M$, then the $Q$-curvature
transforms in the following way
\begin{gather}\label{eq:confP}
   P_g w + 2 Q_g = 2 Q_{\tilde{g}} e^{4 w},
\end{gather}
where $P_g$ is the {\em Paneitz operator}
\[
P_g (\varphi) = \Delta_g^2 \varphi + {\rm div}_g \left(\frac 23 R_g g -
2 {\rm Ric}_g \right) d \varphi,  \qquad  \varphi \in C^\infty(M),
\]
introduced in \cite{p1}. When we modify conformally $g$ as before,
also $P_g$ obeys a simple transformation law, which is the following
\[
  P_{\tilde{g}} = e^{- 4 w} P_g.
\]
The last formula and \eqref{eq:confP} are indeed analogous to
classical ones which hold on compact surfaces: in fact, if we let
$K_g$ denote the Gauss curvature of a surface $\Sigma$ it is well
known that
\begin{gather*}
    \Delta_{\tilde{g}} = e^{-2w} \Delta_g, \qquad  -
    \Delta_g w + K_g = K_{\tilde{g}} e^{2w}.
\end{gather*}
In addition to these, we have an analogy with the classical
Gauss--Bonnet formula ($\chi(\Sig)$ stands here for the Euler
characteristic of $\Sig$)
\begin{gather*}
\int_\Sig K_g dV_g = 2\pi \chi(\Sig).
\end{gather*}
Indeed, if $W_g$ denotes the Weyl tensor of a four manifold $(M,g)$,
one has
\begin{gather}\label{eq:GB4}
    \int_M \left( Q_g + \frac{|W_g|^2}{8} \right) dV_g = 4 \pi^2
    \chi(M).
\end{gather}
Notice that,  from either \eqref{eq:confP} or \eqref{eq:GB4}, it
follows that the total $Q$-curvature of $M$, which we denote by
\begin{gather*}
    k_P = \int_M Q_g dV_g,
\end{gather*}
is a conformal invariant of $(M,g)$: we refer for example to the
survey \cite{cy99} for more details.

Since for two-dimensional surfaces the integral of the Gauss
curvature characterizes completely the topology, it is reasonable to
expect that conformal invariant objects on four mani\-folds, like
$P_g$ or $k_P$, could at least provide some partial geometric
information. Indeed some results in this spirit were obtained by
Gursky in \cite{g2} (see also \cite{g}). If a manifold of
non-negative Yamabe class $Y(g)$ satisf\/ies also $k_P \geq 0$, then
${\rm ker}\, P_g$ consists only of the constant functions and $P_g \geq 0$,
namely $P_g$ is a non-negative operator. If in addition $Y(g) > 0$,
then the f\/irst Betti number of $M$ vanishes, unless $(M,g)$ is
conformally equivalent to a quotient of $S^3 \times \R$. On the
other hand, if $Y(g) \geq 0$, then $k_P \leq 8 \pi^2$ with equality
holding if and only if $(M,g)$ is conformally equivalent to the
standard sphere.

In two dimensions, the {\em uniformization theorem} asserts that
every compact surface carries a conformal metric with constant Gauss
curvature. Such deformation has a strong geometric relevance, since
the deformed metric will be a model one: spherical, Euclidean or
hyperbolic. In four dimensions one might ask the same question for
the $Q$-curvature: this being a scalar quantity, it will not be
suf\/f\/icient to control the whole curvature tensor, however the new
metric will be a privileged one and might have special properties.

Writing $\tilde{g} = e^{2 w} g$, by \eqref{eq:confP} the question
amounts to f\/inding a solution of the equation
\begin{gather}\label{eq:Qc}
    P_g w + 2 Q_g = 2 \ov{Q} e^{4 w},
\end{gather}
where $\ov{Q}$ is a real constant. Problem \eqref{eq:Qc} is
variational, and solutions can be found as critical points of the
following functional
\begin{gather*}
    II(u) = \langle P_g u, u \rangle + 4 \int_M Q_g u dV_g - k_P
    \log \int_M e^{4u} dV_g, \qquad u \in H^{2}(M).
\end{gather*}
Here $H^2(M)$ is the space of real functions on $M$ which are of
class $L^2$ together with their f\/irst and second derivatives, and
the symbol $\langle P_g u, v \rangle$ stands for
\begin{gather*}
    \langle P_g u, v \rangle = \int_M\! \left( \D_g u \D_g v + \frac 23
    R_g \n_g u \cdot \n_g v - 2 ({\rm Ric}_g \n_g u, \n_g v) \right) dV_g,
    \qquad  u, v \in H^2(M).
\end{gather*}
Notice that the functional $II$ is well def\/ined on $H^2(M)$: in
fact, by the Adams inequality (see~\cite{ada, cy95}) one has
\begin{gather}\label{eq:adams}
    \log \int_M e^{4(u - \ov{u})} dV_g \leq \frac{1}{8 \pi^2}
    \langle P_g u, u \rangle + C, \qquad  P_g \geq 0, \qquad
     u \in H^2(M),
\end{gather}
where $\ov{u}$ is the average of $u$ on $M$ and where $C$ depends
only on $M$, so the last term in $II$ is well def\/ined in this
Sobolev space.

Problem \eqref{eq:Qc} has been solved in \cite{cy95} for the case in
which $P_g$ is a non-negative operator and when $k_P < 8 \pi^2$: by
the above-mentioned result of Gursky, suf\/f\/icient conditions for
these assumptions to hold are that $Y(g) \geq 0$ and that $k_P \geq
0$ (and $(M,g)$ is not conformal to the standard sphere).  More
general suf\/f\/icient conditions for the above hypotheses to hold have
been obtained by Gursky and Viaclovsky in \cite{gv}.

Under the assumptions in \cite{cy95}, by \eqref{eq:adams} the
functional $II$ is bounded from below and coercive, hence solutions
can be found as global minima. The result in~\cite{cy95} has also
been extended in~\cite{b1} to higher-dimensional manifolds
(regarding higher-order operators and curvatures) using a geometric
f\/low.

The solvability of \eqref{eq:Qc}, under the above hypotheses, is
useful in the study of some conformally invariant fully non-linear
equations, as shown in \cite{cgy}. Some remarkable geometric
consequences of this study, given in \cite{cgyann,cgy}, are
the following. If a manifold of positive Yamabe class satisf\/ies $k_P
> 0$, then there exists a conformal metric with positive Ricci
tensor, and hence $M$ has f\/inite fundamental group. Furthermore,
under the additional quantitative assumption $k_P > \frac 18 \int_M
|W_g|^2 dV_g$, $M$ must be dif\/feomorphic to the standard four-sphere
or to the standard projective space.

 As we already mentioned, $Q$-curvature (with its
higher-dimensional versions and the Paneitz operator) is relevant in
the study of many other questions, like log-determinant formulas,
compactif\/ication of locally conformally f\/lat manifolds,
Poincar\'e--Einstein metrics, ambient metrics, tractor calculus,
volume renormalization, scattering theory and others, see the
(incomplete) list of references \cite{br, bg, bo, bcy, cqy, cqy2,
fg, fh, go, gope, gjms, graju, gz}.

This note concerns some recent progress about problem \eqref{eq:Qc},
in particular an extension of the {\em uniformization} result of
\cite{cy95}, which is given in~\cite{dm}.

\begin{theorem}\label{th:ex}
Suppose ${\rm ker} \, P_g = \{{\rm const}\}$, and assume that $k_P \neq 8 k
\pi^2$ for $k = 1, 2, \dots$. Then $(M,g)$ admits a conformal metric
with constant $Q$-curvature.
\end{theorem}

The assumptions in Theorem \ref{th:ex} are conformally invariant and
generic, and hence apply to a large class of four manifolds,
especially to some with negative curvature. For example, pro\-ducts of
two negatively-curved surfaces might have total $Q$-curvature
greater than $8 \pi^2$, see \cite{dm1}. However the theorem does not
cover all manifolds, for example conformally f\/lat manifolds with
positive and even Euler characteristic, by \eqref{eq:GB4}, are
excluded.

Our assumptions include those made in \cite{cy95} and one (or both)
of the following two possibilities
\begin{gather}\label{eq:kp3}
    k_P \in (8 k \pi^2, 8 (k+1) \pi^2) \quad \hbox{ for some } k \in \N;
\\
\label{eq:kp4}
    P_g \hbox{ possesses $\ov{k}$ (counted with
    multiplicity)
    negative eigenvalues}.
\end{gather}
In these cases the functional $II$ is unbounded from below, and
hence it is necessary to f\/ind critical points which are possibly of
saddle type. This is done using a new min-max scheme, which we are
going to describe brief\/ly in the next section, depending on $k_P$
and the spectrum of $P_g$ (in particular on the number of negative
eigenvalues $\ov{k}$, counted with multiplicity).

One fundamental issue in applying variational techniques is
compactness: it is possible to prove that indeed solutions of
\eqref{eq:Qc} (or of some perturbation) stay compact provided $P_g$
has no kernel and $k_P$ stay bounded away from $8 \pi^2 \N$. The
following result was proved in \cite{dr, mal} using blow-up
analysis.

\begin{theorem}[\cite{mal}]\label{th:bd}
Suppose ${\rm ker}\, P_g = \{{\rm const}\}$ and that $(u_l)_l$ is a
sequence of solutions to
\begin{gather*}
    P_g u_l + 2 Q_l = 2 k_l e^{4 u_l} \qquad \hbox{in } M,
\end{gather*}
satisfying $\int_M e^{4 u_l} dV_g = 1$, where $k_l = \int_M Q_l d
V_g$, and where $Q_l \to Q_0$ in $C^0(M)$. Assume also that
\begin{gather*}
    k_0 := \int_M Q_0 dV_g \neq 8 k \pi^2 \qquad
    \hbox{ for } k = 1, 2, \dots.
\end{gather*}
Then $(u_l)_l$ is bounded in $C^\a(M)$ for any $\a \in (0,1)$.
\end{theorem}

We will show below how this result is useful in deriving existence.
For reasons of brevity, we will describe the proof of Theorem
\ref{th:ex} only: as mentioned, we will f\/ind in general solutions as
saddle points of $II$. To use variational techniques we will try to
understand the topological properties of the sublevels, in
particular the very negative ones.

If $k_P > 8 \pi^2$, it is easy to see that $II$ is unbounded from
below: one can f\/ix a point $x \in M$ and consider a test function
with the following form
\[
  \var_{\l,x}(y) \simeq \frac 14 \log \left( \frac{2 \l}{1 + \l^2 {\rm dist}(y,x)^2}
  \right), \qquad y \in M,
\]
where ${\rm dist}(\cdot,\cdot)$ denotes the metric distance on $M$. For
$\l \to + \infty$ we have that $II(\var_{\l,x}) \to - \infty$.
Another way to attain large negative value of $II$ is to consider
negative eigenvectors of $P_g$ with large norm: in the next section
we will focus on the f\/irst alternative only, and reason looking at
how the function $e^{4u}$ is {\em distributed} over dif\/ferent
regions of $M$.

The last section describes other  related results and some open
problems in this direction.

\section{Sketch of the proof}\label{s:proof}

We give here the main ideas for the proof of Theorem \ref{th:ex}:
throughout this section we assume for simplicity that $P_g$ is
positive def\/inite (except on constants), referring to \cite{dm} for
details when negative eigenvalues are present.

\subsection{Improved Adams inequality and applications}\label{ss:ia}

One of the main ingredients in our proof is an improvement of the
Adams inequality on the functions $u$ for which $e^{4u}$ is {\em
spread} into multiple regions: here is the precise statement.

\begin{lemma}\label{l:imprc}
For a fixed integer $\ell$, let $\Omega_1, \dots, \Omega_{\ell+1}$
be subsets of $M$ satisfying ${\rm dist}(\Omega_i,\Omega_j) \geq \d _0$
for $i \neq j$, where $\d_0$ is a positive real number, and let
$\g_0 \in \left( 0, \frac{1}{\ell+1} \right)$. Then, for any
$\tilde{\e}> 0$ there exists a constant $C = C(\ell, \tilde{\e},
\d_0, \g_0)$ such that
\[
  \log \int_M e^{4(u - \ov{u})} dV_g \leq C + \frac{1}{8 (\ell+1)
  \pi^2 - \tilde{\e}} \langle P_g u, u \rangle
\]
for all the functions $u \in H^2(M)$ satisfying
\begin{gather}\label{eq:ddmm}
    \frac{\int_{\Omega_i} e^{4u} dV_g}{\int_M e^{4u} dV_g} \geq \g_0,
  \qquad \forall \; i \in \{1, \dots, \ell+1\}.
\end{gather}
\end{lemma}

The original improvement argument is given in
\cite{cl}, where the authors treat $H^1$ functions on surfaces and
the case $\ell = 1$: some modif\/ications are needed here to deal with
the fourth order operator. The main step of the proof consists in
constructing cutof\/f functions $g_i$ which are identically equal to
$1$ on $\Omega_i$ and which have mutually disjoint supports. Then
one applies \eqref{eq:adams} to $g_i u$ and chooses the index $i$ so
that $\int_M g_i^2 (\D u)^2 dV_g$ is minimal: the remaining terms of
$\int_M (\D (g_i u))^2 dV_g$ can be treated as perturbations and
yield $o(\|u\|_{H^2(M)})$.

The above lemma states that if $k_P < 8 (k+1) \pi^2$ for some
natural $k$, and if $u$ satisf\/ies \eqref{eq:ddmm}  with $\ell = k$, then $II(u)$
stays bounded from below. Since we can choose $\d_0$ and the
$\Omega_i$'s arbitrarily, this suggests that if $II(u)$ is
suf\/f\/iciently low, then $e^{4 u}$ has to {\em concentrate} near at
most $k$ points of~$M$. Using a covering argument this can indeed be
made rigorous, yielding the following result.

\begin{lemma}\label{l:II<-M}
Assuming $P_g \geq 0$ and $k_P \in (8 k \pi^2, 8(k+1) \pi^2)$ with
$k \geq 1$, the following property holds. For any  $\e
> 0$ and any $r > 0$ there exists a large positive $L = L(\e, r)$ such that for
every $u \in H^2(M)$ with $II(u) \leq - L$ there exist $k$ points
$p_{1,u}, \dots, p_{k,u} \in M$ such that
\begin{gather*}
 \int_{M \setminus \cup_{i=1}^{k} B_r(p_{i,u})} e^{4u} dV_g < \e.
\end{gather*}
\end{lemma}

\begin{figure}[t]
\centering
 \includegraphics[width=6cm]{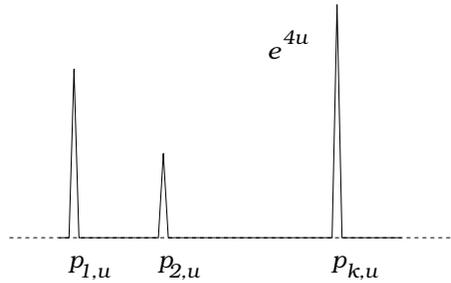}
\caption{Shape of functions with low energy.}\label{fig:1}
\end{figure}

Lemma \ref{l:II<-M} basically states that if $II(u)$ is large
negative, then the conformal volume $e^{4u}$ has a shape as in
Fig.~\ref{fig:1}.

If we assume the normalization $\int_M e^{4u} dV_g = 1$, the f\/igure
suggests to consider the following set
\[
  M_k = \left\{ \sum_{i=1}^k t_i \d_{x_i}  : \; x_i \in M, \; t_i \in [0,1],
  \sum_{i=1}^k t_i = 1 \right\},
\]
known in literature as the set of {\em formal barycenters} of $M$.
For $k = 1$ $M_1$ is homeomorphic to $M$, but for larger $k$'s this
is a {\em stratified set}, namely union of open manifolds of
dif\/ferent dimensions, whose maximal one is $5k -1$. For further
properties of these sets we refer the reader to~\cite{bc,bred}: we are going to need especially the following one,
which can be proved with an argument in algebraic topology.

\begin{lemma}[well-known]\label{l:nonc}
If $M$ is a compact manifold, for any $k \geq 1$ the set $M_k$ is
non-contractible.
\end{lemma}

Using Lemma \ref{l:II<-M} one can construct a continuous (and
non-trivial) map from very negative sublevels of $II$ into $M_k$.

\begin{proposition}\label{p:map}
For $k \geq 1$ (see \eqref{eq:kp3}) there exists a large $L > 0$ and
a continuous map $\Psi$ from the sublevel $\{ II \leq - L \}$ into
$M_k$ which is topologically non-trivial.
\end{proposition}

\begin{figure}[t]
\centering
\includegraphics[width=10cm]{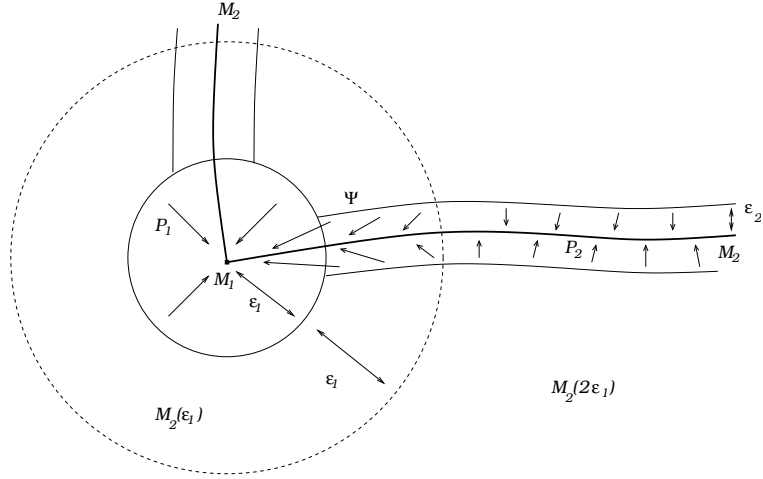}
\caption{A sketch on the construction of ${\Psi}$ for $k = 2$.}\label{fig:2}
\end{figure}

By {\em topologically non-trivial} we mean that the action of
$\Psi_*$ on the homology classes of the sublevel $\{ II \leq - L \}$
has non-trivial image. The proof of the latter proposition is rather
involved and takes some long technical work: we limit ourselves to
explain some of the ideas for the construction. It is easy to {\em
project} when $e^{4u}$ is {\em close} to smooth pieces of $M_k$, but
some dif\/f\/iculties may occur when approaching the singularities. In
this case, one can project f\/irst on the subsets of lowest possible
dimension, and then take homotopies with the higher-dimensional
ones, see Fig.~\ref{fig:2}.

\subsection{The min-max scheme}

In order to proceed we need some preliminary notation. For $\d
> 0$ small, consider a smooth non-decreasing cut-of\/f function
$\chi_\d : \R_+ \to \R$ satisfying the following properties
\begin{gather*}
 \chi_\d(t) =  t \qquad \hbox{for } t \in [0,\d]; \nonumber\\
\chi_\d(t) =  2 \d \qquad \hbox{for } t \geq 2 \d; \\ 
\chi_\d(t) \in [\d, 2 \d] \qquad \hbox{for } t \in [\d, 2 \d].
\nonumber
\end{gather*}
Then, given $\s \in M_k$ $\left( \s = \sum_{i=1}^k t_i \d_{x_i}
\right)$ and $\l
> 0$, we def\/ine the function $\var_{\l,\s} : M \to \R$ as
\begin{gather}\label{eq:pls}
  \var_{\l,\s} (y) = \frac 14 \log
  \sum_{i=1}^k t_i \left( \frac{2 \l}{1 + \l^2 \chi_\d^2
  \left( d_i(y) \right)} \right)^4, \qquad  y \in M,
\end{gather}
where  $d_i(y) = {\rm dist}(y,x_i)$, $y \in M$, with ${\rm dist}(\cdot, \cdot)$
denoting the distance function on $M$. We have then the following
result.

\begin{proposition}\label{p:sublev}
Let $\var_{\l,\s}$ be defined as in \eqref{eq:pls}. Then, as $\l \to
+ \infty$ the following properties hold
\begin{enumerate}\itemsep=0pt

\item[(i)] $e^{4 \var_{\l,\s}} \rightharpoonup \s$ weakly in the sense
of distributions;

\item[(ii)] $II(\var_{\l,\s}) \to - \infty$ uniformly for $\s
\in M_k$.
\end{enumerate}
Moreover, if $\Psi$ is given in Proposition {\rm \ref{p:map}}, the map
$\sigma \mapsto \var_{\l,\s} \mapsto \Psi(\var_{\l,\s})$ converges
to the identity on $M_k$ as $\l$ tends to infinity.
\end{proposition}

It is indeed the last property, together with the
non-contractibility of $M_k$,  which guarantees the non-triviality
of the map $\Psi$.

We next introduce the variational scheme which provides existence of
solutions for~\eqref{eq:Qc}. Let $\widehat{M_k}$ denote the
(contractible) topological cone over $M_k$, which can be represented
as $\widehat{M_k} = M_k \times [0,1]$ with $M_k \times \{0\}$
collapsed to a single point. Let f\/irst $L$ be so large that
Proposition \ref{p:map} applies with $\frac L4$, and then let
$\ov{\l}$ be so large that Proposition \ref{p:sublev} applies for
this value of $L$. Fixing this number $\ov{\l}$, we def\/ine the
following class of maps
\begin{gather*}
    \Pi_{\ov{\l}} = \left\{ \pi : \widehat{M_k} \to H^2(M)
  : \; \pi \hbox{ is continuous and } \pi(\cdot \times \{1\})
 = \var_{\ov{\l},\cdot} \hbox{ on } M_k \right\}.
\end{gather*}
We have then the following properties.

\begin{lemma}\label{l:min-max}
The set $\Pi_{\ov{\l}}$  is non-empty and moreover, letting
\begin{gather}\label{eq:gapflow}
    \ov{\Pi}_{\ov{\l}} = \inf_{\pi \in \Pi_{\ov{\l}}}
  \; \sup_{m \in \widehat{M_k}} II(\pi(m)), \qquad
  \hbox{ one has } \qquad \ov{\Pi}_{\ov{\l}} > - \frac
  L2.
\end{gather}
As a consequence, the functional $II$ possesses a Palais--Smale
sequence at level $\ov{\Pi}_{\ov{\l}}$.
\end{lemma}

To check that $\Pi_{\ov{\l}}$ is non empty, it is suf\/f\/icient to
consider the map $\pi(t, \s) = t \var_{\ov{\l},\s}$, which is well
def\/ined on $\widehat{M_k}$. A {\em Palais--Smale sequence at level}
$c$ is a sequence $(u_n)_n$ such that $II(u_n) \to c$ and such that
$II'(u_n) \to 0$ as $n \to + \infty$. To produce it one considers
any map in $\Pi_{\ov{\l}}$ with $\sup_{m \in \widehat{M_k}}
II(\pi(m))$ suf\/f\/iciently close to $\ov{\Pi}_{\ov{\l}}$ and let it
evolve via the gradient f\/low of $II$, keeping its boundary f\/ixed. By
\eqref{eq:gapflow} the maximum value of $II$ on the evolved set,
which still belongs to $\Pi_{\ov{\l}}$, cannot go below
$\ov{\Pi}_{\ov{\l}}$, so the gradient of $II$ along the evolution
has to become small somewhere along the f\/low, see for example
\cite{strb}.

If Palais--Smale sequences are bounded, then it is easy to see that
they have to converge to a~solution of \eqref{eq:Qc}. Unfortunately
boundedness is not known in general, so we tackle this problem using
an argument due to M.~Struwe, see~\cite{str}. For $\rho$ in a
neighborhood of $1$ we def\/ine the functional $II_{\rho} : H^2(M) \to
\R$ by
\begin{gather*}
  II_\rho(u) = \langle P_g u, u \rangle + 4 \rho \int_M Q_g d V_g
  - 4 \rho k_P \log \int_M e^{4 u} d V_g, \qquad u \in H^2(M),
\end{gather*}
whose critical points give rise to solutions of the equation
\begin{gather}\label{eq:mod}
    P_g u + 2 \rho Q_g = 2 \rho k_P e^{4 u} \qquad \hbox{in } M.
\end{gather}
Running the above variational scheme one f\/inds a min-max value as in
Lemma~\ref{l:min-max}, which we call~$\ov{\Pi}_\rho$, and a
corresponding Palais--Smale sequence. One can easily show that
$\frac{\ov{\Pi}_\rho}{\rho}$ is non-increasing, and hence
dif\/ferentiable, in $[1 - \rho_0, 1 + \rho_0]$ provided $\rho_0$ is
suf\/f\/iciently small. Using the following result, which can be proved
as in \cite{djlw}, by the above comments and Theorem \ref{th:bd} one
deduces then existence of solutions to~\eqref{eq:Qc}.

\begin{lemma}\label{c:c}
Let $\L \subset [1 - \rho_0, 1 + \rho_0]$ be the (dense) set of
$\rho$ for which the function $\frac{\ov{\Pi}_\rho}{\rho}$ is
differen\-tiable. Then for $\rho \in \L$  $II_\rho$ possesses a
bounded Palais--Smale sequence $(u_l)_l$ at level $\ov{\Pi}_{\rho}$.
\end{lemma}

A particular case of the above variational method has been used in~\cite{djlw} to study a mean f\/ield equations on compact surfaces,
which in our setting would correspond to the case $k = 1$.

\section{Final remarks}

The set of barycenters of a manifold (or of a set)  has been used
crucially in literature for the study of problems with lack of
compactness, see \cite{bah,bc}. In particular, for
Yamabe-type equations (including the Yamabe equation and several
other applications), it has been employed to understand the
structure of the {\em critical points at infinity} (or asymptotes)
of the Euler functional, namely the way compactness is lost through
a pseudo-gradient f\/low. Our use of the set~$M_k$, although the map~$\Psi$  presents some analogies with the Yamabe case, is of
dif\/ferent type since it is employed to reach low energy levels and
not to study critical points at inf\/inity. As mentioned above, we
consider a projection onto the $k$-barycenters $M_k$, but starting
only from functions in $\{ II \leq - L \}$, whose concentration
behavior is not as clear as that of the asymptotes for the Yamabe
equation.

Theorem \ref{th:ex} has been recently extended to all dimensions in
\cite{nd1}: one of the issues there is that the explicit expression
of the Paneitz operator is not known in general, but only the
principal part is needed. Also, the blow-up analysis to obtain the
counterpart of Theorem \ref{th:bd} relies almost entirely on the
Green's representation formula.

The above approach has been also used in \cite{nd2, nd3} to study
the case of four manifolds with boundary: in fact, in two papers by
A.~Chang and J.~Qing, see \cite{cq1,cq2}, it was constructed a
scalar quantity, the $T$-{\em curvature}, which represents an
analogue of the geodesic curvature for surfaces with boundary. If
one wants to prescribe constant $Q$-curvature in the interior and
zero $T$-curvature (and mean curvature) on the boundary, it turns
out that the quantization of conformal volume at the boundary is one
half of the one in the interior. As a consequence, one needs to
modify the topological argument and consider formal barycenters
which allow a double number of Dirac deltas on $\partial M$ compared
to the interior.

We also mention the papers \cite{b2, ms, wx}, where the problem of
prescribing non-constant $Q$-curvature on spheres was considered. In
this case the analysis is more delicate since the problem is non
compact, due to the presence of the M\"obius group. One has then to
analyze how compactness is lost and employ topological arguments,
like degree of Morse theory, to derive existence results. In fact
Morse theory can be used, jointly with the Poincar\'e--Hopf theorem,
to compute the Leray--Schauder degree of \eqref{eq:Qc}. By Theorem
\ref{th:bd}, if $P_g$ has no kernel and if $k_P \not\in 8 \pi^2 \N$,
then the solutions of \eqref{eq:Qc} stay uniformly bounded, so the
degree of the equation is well def\/ined, and in \cite{maldeg} it was
proven that it is given by
\begin{gather*}
        (-1)^{\ov{k}} \qquad  \hbox{for } k_P < 8 \pi^2; \\
        (-1)^{\ov{k}} (-\chi(M)+1) \cdots (-\chi(M)+k)/k!
        \qquad
       \hbox{for } k_P \in (8 k \pi^2, 8 (k+1) \pi^2), \quad k \in \N,
\end{gather*}
where $\ov{k}$ stands for the number of negative eigenvalues of
$P_g$, counted with multiplicity. This result is related to one
obtained in \cite{cl2} for a mean f\/ield equation on compact
surfaces, and based on sophisticated blow-up techniques: our proof
is of dif\/ferent type and relies instead on Morse theory, as
mentioned.

A natural question one may ask is whether Theorem \ref{th:ex} can be
extended to the situations in which $P_g$ has non-trivial kernel or
when $k_P$ is an integer multiple of~$8 \pi^2$: in this cases more
precise characterizations of $P_g$ and $Q_g$ (or on their mutual
relation) should be given.

One simple situation to consider is when $P_g$ has non-trivial
kernel and $k_P = 0$. In this case~\eqref{eq:Qc} has the following
simple expression
\begin{gather*}
   P_g u + 2 Q_g = 0.
\end{gather*}
Clearly, by Fredholm's alternative, the equation is solvable if and
only if $Q_g$ is perpendicular to ${\rm ker}\, P_g$. It seems possible to
violate this condition but no example so far is known.

If one wishes to study the case of non-trivial kernel in more
generality, possibly employing min-max arguments, it would be
necessary to understand the behavior of $II$ on ${\rm ker}\, P_g$, in
particular on functions with large norm. For example, taking $v \in
{\rm ker}\, P_g$, $v \neq 0$, one has
\[
  II(t v) = 4 t \int_M Q v d V_g - k_P \log \int_M e^{4tv} dV_g.
\]
Since every element in ${\rm ker}\, P_g$ is smooth we have $v \in
L^\infty(M)$, and in the limit $t \to + \infty$ one has $\log \int_M
e^{4tv} dV_g \simeq 4 t \max_M v$. It follows that for $v \in {\rm ker}\,
P_g$
\[
  II(t v) = 4 t \left(  \int_M Q v d V_g - k_P \max_M v \right) + o(t),
  \qquad t \to + \infty.
\]
It would be useful to prove that this function either has constant
sign on the unit sphere of ${\rm ker}\, P_g$, or that its negative sublevel
has non-trivial topology. Under these conditions one might hope to
adapt the variational scheme in \cite{dm} and prove existence of
solutions to~\eqref{eq:Qc}.

Also the case $k_p = 8 k \pi^2$ should require some new geometric
insight. In the recent paper~\cite{llp} the authors considered four
manifolds such that $k_p = 8 \pi^2$. They approached the problem
perturbing the equation as in \eqref{eq:mod}, letting $\rho$ tending
to $1$ from below: indeed for $\rho < 1$ by the Adams inequality
$II_\rho$ admits a minimizer $u_\rho$. In \cite{llp} the asymptotic
behavior of $u_\rho$ was studied and it was shown that if the
function
\begin{gather}\label{eq:funct}
    (\D_{g,y} S(x,y))_{y=x} + 4 |\n_{g,y} S(x,y)_{y=x}|^2 - \frac{R_g(x)}{8}
\end{gather}
is positive on $M$, then $u_\rho$ converges to a solution of
\eqref{eq:Qc}. In the above formula $S$ stands for the regular part
of the Green's function of $P_g$: precisely, using the convention
that
\[
  P_g G(x,y) + 2 Q_g(y) = 16 \pi^2 \d_x,
\]
one has that, near the point $x$, $G$ satisf\/ies
\[
  G(x,y) = - 2 \log {\rm dist}(x,y) + S(x,y)
\]
for some $C^2$-smooth $S$. However, there are so far no known
conditions which can guarantee the positivity of the function in
\eqref{eq:funct}. A result of this kind would be a counterpart of
the {\em positive mass theorem} by Schoen and Yau for compact
manifolds of positive scalar curvature.

\subsection*{Acknowledgements}

The author has been supported by M.U.R.S.T within the PRIN 2006 {\em
Variational Methods and Nonlinear Differential Equations}.

\pdfbookmark[1]{References}{ref}
\LastPageEnding

\end{document}